\documentclass{amsart}
\usepackage{amsfonts}
\usepackage{amssymb}
\usepackage{amscd}
\usepackage{graphicx}
\usepackage{xypic}
\usepackage{psfrag}

\psfrag{rij-1-rij}{$r_{i,j-1} - r_{ij}$}
\psfrag{ri+1j-rij}{$r_{i+1,j} - r_{ij}$}
\psfrag{R01}{$R_{01}$}
\psfrag{R12}{$R_{12}$}
\psfrag{Rn-1n}{$R_{n-1,n}$}
\psfrag{sig1}{$\sigma_1$}
\psfrag{sig2}{$\sigma_2$}
\psfrag{tau1}{$\tau_1$}
\psfrag{taun-1}{$\tau_{n-1}$}
\psfrag{Ri-1i}{$R_{i-1,i}$}
\psfrag{sigi}{$\sigma_i$}
\psfrag{taui-1}{$\tau_{i-1}$}
\psfrag{U}{$U$}
\psfrag{D}{$D$}
\psfrag{R}{$R$}
\psfrag{P}{$P$}
\psfrag{Q0}{$Q_0$}
\psfrag{P0}{$P_0$}
\psfrag{Z}{$Z$}
\psfrag{cdot}{$\cdot$}
\psfrag{=}{$=$}
\psfrag{Ptil}{$\Tilde P$}
\psfrag{sigma}{$\sigma$}
\psfrag{tau}{$\tau$}
\psfrag{mu}{$\mu$}
\psfrag{R}{$R$}
\psfrag{T}{$T$}
\psfrag{Ti-1i}{$T_{i-1,i}$}
\psfrag{L}{$L$}
\psfrag{Q}{$Q$}
\psfrag{Q=}{$Q =$}
\psfrag{=P}{$= P$}
\psfrag{Qtil}{$\Tilde Q$}
\psfrag{W}{$W$}
\psfrag{W0}{$W_0$}
\psfrag{Wtil}{$\Tilde W$}
\psfrag{A}{$A$}
\psfrag{B*C}{$B \cdot C$}
\psfrag{A*B*C}{$A \cdot B \cdot C$}
\psfrag{X}{$X$}
\psfrag{Y}{$Y$}
\psfrag{S}{$S$}
\psfrag{V}{$V$}
\psfrag{W1}{$W_1$}
\psfrag{W2}{$W_2$}
\psfrag{W3}{$W_3$}
\psfrag{Wn}{$W_n$}
\psfrag{P1}{$P_1$}
\psfrag{P2}{$P_2$}
\psfrag{P3}{$P_3$}
\psfrag{Pn}{$P_n$}
\psfrag{Q1}{$Q_1$}
\psfrag{Q2}{$Q_2$}
\psfrag{Qn-1}{$Q_{n-1}$}
\psfrag{P1*Q1}{$P_1 \cdot Q_1$}
\psfrag{P2*Q2}{$P_2 \cdot Q_2$}
\psfrag{B}{$B$}
\psfrag{E}{$E$}
\psfrag{F}{$F$}
\psfrag{x}{$x$}
\psfrag{c1}{$c_1$}
\psfrag{c2}{$c_2$}
\psfrag{c3}{$c_3$}
\psfrag{c4}{$c_4$}
\psfrag{ck}{$c_k$}
\psfrag{d1}{$d_1$}
\psfrag{d2}{$d_2$}
\psfrag{dl}{$d_l$}
\psfrag{C}{$C$}
\psfrag{gamma}{$\gamma$}
\psfrag{gamma'}{$\gamma'$}
\psfrag{X'}{$X'$}
\psfrag{Y'}{$Y'$}
\psfrag{X=}{$X =$}
\psfrag{X0}{$X_0$}
\psfrag{Xtil}{$\Tilde X$}
\psfrag{=Y}{$= Y$}
\psfrag{Y0}{$Y_0$}
\psfrag{Ytil}{$\Tilde Y$}
\psfrag{A0}{$A_0$}
\psfrag{Atil}{$\Tilde A$}
\psfrag{T'}{$T'$}

\newcommand{\Z}{{\mathbb Z}}

\newcommand{\PP}{{\mathcal P}}

\newcommand{\Th}{{\text{th}}}
\newcommand{\St}{{\text{st}}}

\newcommand{\bull}{{\sssize \bullet}}
\newcommand{\rank}{\operatorname{rank}}

\newcommand{\rect}{\operatorname{rect}}

\newcommand{\picC}[1]{\includegraphics[scale=0.50]{#1.eps}}

\newcommand{\picT}[1]{\includegraphics[scale=0.60]{#1.eps}}

\newcommand{\pathskip}{\vspace{5pt}}

\newcommand{\attach}[3]{\frac{#1|#3}{#2}}

\newcommand{\PART}[1]{\includegraphics[scale=.50]{p#1.eps}}

\newcommand{\TABLr}[2]{\raisebox{#2pt}{\includegraphics[scale=.60]{t#1.eps}}}
\newcommand{\TABL}[1]{\includegraphics[scale=.60]{t#1.eps}}

\theoremstyle{plain}

\newtheorem{thm}{Theorem}
\newtheorem{lemma}{Lemma}
\newtheorem{prop}{Proposition}
\newtheorem*{cor}{Corollary}
\newtheorem{conj}{Conjecture}
\newtheorem*{conj1a}{Conjecture 1A}

\theoremstyle{definition}
\newtheorem{example}{Example}

\newcommand{\refthm}[1]{Theorem~\ref{#1}}
\newcommand{\refprop}[1]{Proposition~\ref{#1}}
\newcommand{\reflemma}[1]{Lemma~\ref{#1}}

\newcommand{\refconj}[1]{Conjecture~\ref{#1}}

\newenvironment{mathenum}{\begin{enumerate}}{\end{enumerate}}

\newenvironment{romenum}{\begin{enumerate}}{\end{enumerate}}

\title{On a conjectured formula for quiver varieties}
\date{\today}
\author{Anders Skovsted Buch}
\address{Massachusetts Institute of Technology \\
  Building 2, Room 248 \\
  77 Massachusetts Avenue \\
  Cambridge, MA 02139
}
\email{abuch@math.mit.edu} 

\begin{document}
\maketitle

\section{Introduction}

The goal of this paper is to prove some combinatorial results about a
formula for quiver varieties given in \cite{buch.fulton:chern}.

Let $X$ be a non-singular complex variety and $E_0 \to E_1 \to E_2 \to
\cdots \to E_n$ a sequence of vector bundles and bundle maps over $X$.
A set of {\em rank conditions\/} for this sequence is a collection of
non-negative integers $r = (r_{ij})$ for $0 \leq i < j \leq n$.  This
data defines a degeneracy locus in $X$,
\[ \Omega_r(E_\bull) = \{ x \in X \mid \rank(E_i(x) \to E_j(x)) \leq
    r_{ij} ~\forall i < j \} \,. 
\]

Let $r_{ii}$ denote the rank of the bundle $E_i$.  We will demand that
the rank conditions can {\em occur\/}, i.e.\ that there exists a
sequence of vector spaces and linear maps $V_0 \to V_1 \to \dots \to
V_n$ so that $\dim(V_i) = r_{ii}$ and $\rank(V_i \to V_j) = r_{ij}$.
This is equivalent to the conditions $r_{ij} \leq \min(r_{i,j-1},
r_{i+1,j})$ for $i < j$, and $r_{ij} - r_{i,j-1} - r_{i+1,j} +
r_{i+1,j-1} \geq 0$ for $j-i \geq 2$.

Given two vector bundles $E$ and $F$ on $X$ and a partition $\lambda$,
we let $s_\lambda(F - E)$ denote the super-symmetric Schur polynomial
in the Chern roots of these bundles.  By definition this is the
determinant of the matrix whose $(i,j)^\Th$ entry is the coefficient
of the term of degree $\lambda_i + j - i$ in the formal power series
expansion of the quotient of total Chern polynomials
$c_t(E^\vee)/c_t(F^\vee)$.

The expected (and maximal) codimension for the locus
$\Omega_r(E_\bull)$ in $X$ is
\[ d(r) = \sum_{i<j} (r_{i,j-1} - r_{ij}) \cdot (r_{i+1,j} - r_{ij}) \,. 
\]
The main result of \cite{buch.fulton:chern} gives a formula for the
cohomology class of $\Omega_r(E_\bull)$ when it has this codimension:
\[ [\Omega_r(E_\bull)] = \sum_\mu c_\mu(r)\, s_{\mu_1}(E_1 - E_0) \cdots
   s_{\mu_n}(E_n - E_{n-1}) \,.
\]
Here the sum is over sequences of partitions $\mu = (\mu_1, \dots,
\mu_n)$; the coefficients $c_\mu(r)$ are certain integers given by an
explicit combinatorial algorithm which is described in Section
\ref{sec_algo}.  These coefficients are known to generalize
Littlewood-Richardson coefficients as well as the coefficients
appearing in Stanley symmetric functions \cite{buch.fulton:chern},
\cite{buch:stanley}.  The formula specializes to give new expressions
for all known types of Schubert polynomials \cite{fulton:universal}.

There is no immediate geometric reason for the products of Schur
polynomials appearing in the formula.  However, it is even more
surprising that the coefficients $c_\mu(r)$ all seem to be
non-negative.  Attempts to prove this has led to a conjecture saying
that these coefficients count the number of different sequences of
tableaux satisfying certain conditions \cite{buch.fulton:chern}.
These sequences are called {\em factor sequences\/} and are defined in
Section \ref{sec_algo}.

The main result in this paper is a proof of this conjecture in some
special cases which include all situations where the sequence
$E_\bull$ has up to four bundles.  We will also show that the
conjecture follows from a stronger but simpler conjecture, for which
substantial computational verification has been obtained.  For both of
these results, a sign-reversing involution on pairs of tableaux
constructed by S.~Fomin plays a fundamental role.

In Section \ref{sec_algo} we will explain the algorithm for computing
the coefficients $c_\mu(r)$, as well as the conjectured formula for
these coefficients.  In Section \ref{sec_criterion} we will prove a
useful criterion for recognizing factor sequences.  Section
\ref{sec_fomin} gives an account of Fomin's involution, which in
Section \ref{sec_stronger} is used to formulate the stronger
conjecture mentioned above.  Finally, Section \ref{sec_fourbdl}
contains a proof of this stronger conjecture in special cases.

The work described in this paper can be viewed as a continuation of a
joint geometric project with W.~Fulton, which resulted in the quiver
formula described in \cite{buch.fulton:chern}.  We would like to thank
him for introducing us to the subject of degeneracy loci during this
very pleasant collaboration, and also for numerous suggestions, ideas,
comments, etc.\ during the work on this paper.  We are also extremely
grateful to S.~Fomin who provided the vital involution mentioned
above, and who also collaborated with us in the attempts to prove the
conjecture.


\section{Description of the algorithm}
\label{sec_algo}

This section explains the algorithm for computing the coefficients
$c_\mu(r)$ as well as the conjecture for these coefficients.  We will
first explain this in the ordinary case described in the introduction.
Then we will extend the notions to a more general situation, which for
many purposes is easier to work with.

We will need some notation.  Let $\Lambda = \Z[h_1, h_2, \dots]$ be
the ring of symmetric functions.  The variable $h_i$ may be identified
with the complete symmetric function of degree $i$.  If $I = (a_1,
a_2, \dots, a_p)$ is a sequence of integers, define the Schur function
$s_I \in \Lambda$ to be the determinant of the $p \times p$ matrix
whose $(i,j)^\Th$ entry is $h_{a_i + j - i}$:
\[ s_I = \det(h_{a_i + j - i})_{1 \leq i, j \leq p} \,. \]
(Here one sets $h_0 = 1$ and $h_{-q} = 0$ for $q > 0$.)  A Schur
function is always equal to either zero or plus or minus a Schur
function $s_\lambda$ for a partition $\lambda$.  This follows from
interchanging the rows of the matrix defining $s_I$.  Furthermore, the
Schur functions given by partitions form a basis for the ring of
symmetric functions \cite{macdonald:symmetric*1}, \cite{fulton:young}.

We will give the algorithm for computing the coefficients $c_\mu(r)$
by constructing an element $P_r$ in the $n^\Th$ tensor power of the
ring of symmetric functions $\Lambda^{\otimes n}$, such that
\[ P_r = \sum_\mu c_\mu(r) \, 
   s_{\mu_1} \otimes \dots \otimes s_{\mu_n} \,. 
\]

It is convenient to arrange the rank conditions in a {\em rank
diagram\/}:
\[ \begin{matrix}
E_0 & \to & E_1 & \to & E_2 & \to & \cdots & \to & E_n
\vspace{0.1cm} \\
r_{00} && r_{11} && r_{22} && \cdots && r_{nn} \\
& r_{01} && r_{12} && \cdots && r_{n-1,n} \\
&& r_{02} && \cdots && r_{n-2,n} \\
&&& \ddots \\
&&&& r_{0n}
\end{matrix} \]
In this diagram we replace each small triangle of numbers
\[ \begin{matrix}
r_{i,j-1} && r_{i+1,j} \\
& r_{ij}
\end{matrix} \]
by a rectangle $R_{ij}$ with $r_{i+1,j} - r_{ij}$ rows and $r_{i,j-1}
- r_{ij}$ columns.
\[ R_{ij} = \raisebox{-18pt}{\picC{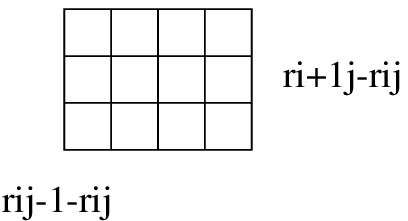}} \]
These rectangles are then arranged in a {\em rectangle diagram\/}:
\[ \begin{matrix}
R_{01} && R_{12} && \cdots && R_{n-1,n} \\
& R_{02} && \cdots && R_{n-2,n} \\
&& \ddots \\
&&& R_{0n}
\end{matrix} \]

It turns out that the information carried by the rank conditions is
very well represented in this diagram.  First, the expected
codimension $d(r)$ for the locus $\Omega_r(E_\bull)$ is equal to the
total number of boxes in the rectangle diagram.  Furthermore, the
condition that the rank conditions can occur is equivalent to saying
that the rectangles get narrower when one travels south-west, while
they get shorter when one travels south-east.  Finally, the element
$P_r$ depends only on the rectangle diagram.

We will define $P_r \in \Lambda^{\otimes n}$ by induction on $n$.
When $n = 1$ (corresponding to a sequence of two vector bundles), the
rectangle diagram has only one rectangle $R = R_{01}$.  In this case
we set
\[ P_r = s_R \in \Lambda^{\otimes 1} \]
where $R$ is identified with the partition for which it is the Young
diagram.  This case recovers the Giambelli-Thom-Porteous formula.

If $n \geq 2$ we let $\Bar r$ denote the bottom $n$ rows of the rank
diagram.  Then $\Bar r$ is a valid set of rank conditions, so by
induction we can assume that 
\begin{equation}
\label{eqn_prbar}
P_{\Bar r} = \sum_\mu c_\mu(\Bar r) \, 
s_{\mu_1} \otimes \dots \otimes s_{\mu_{n-1}}
\end{equation}
is a well defined element of $\Lambda^{\otimes n-1}$.  Now $P_r$ is
obtained from $P_{\Bar r}$ by replacing each basis element $s_{\mu_1}
\otimes \dots \otimes s_{\mu_{n-1}}$ in (\ref{eqn_prbar}) with the sum
\[ \sum_{\stackrel{\sigma_1,\dots,\sigma_{n-1}}{\tau_1,\dots,\tau_{n-1}}}
   \left(\prod_{i=1}^{n-1} c^{\mu_i}_{\sigma_i \tau_i}\right)
   s_{\picC{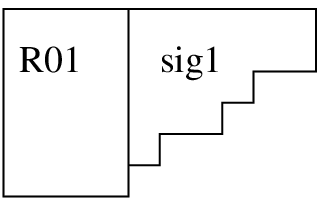}} \otimes \cdots \otimes
   s_{\picC{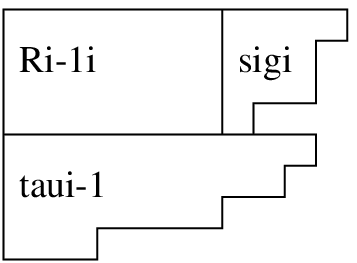}} \otimes \cdots \otimes
   s_{\picC{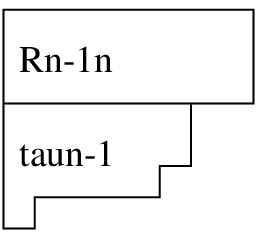}} \,.
\]
This sum is over all partitions $\sigma_1, \dots, \sigma_{n-1}$ and
$\tau_1, \dots, \tau_{n-1}$ such that $\sigma_i$ has fewer rows than
$R_{i-1,i}$ and each Littlewood-Richardson coefficient
$c^{\mu_i}_{\sigma_i \tau_i}$ is non-zero.  A diagram consisting of a
rectangle $R_{i-1,i}$ with (the Young diagram of) a partition
$\sigma_i$ attached to its right side, and $\tau_{i-1}$ attached
beneath should be interpreted as the sequence of integers giving the
number of boxes in each row of this diagram.

It can happen that the rectangle $R_{i-1,i}$ is empty, since the
number of rows or columns can be zero.  If the number of rows is zero,
then $\sigma_i$ is required to be empty, and the diagram is the Young
diagram of $\tau_{i-1}$.  If the number of columns is zero, then the
algorithm requires that the length of $\sigma_i$ is at most equal to
the number of rows $r_{ii} - r_{i-1,i}$ of $R_{i-1,i}$, and the
diagram consists of $\sigma_i$ in the top $r_{ii} - r_{i-1,i}$ rows
and $\tau_{i-1}$ below this, possibly with some zero-length rows in
between.

Next we will describe the conjectured formula for the coefficients
$c_\mu(r)$.  We will need the notions of (semistandard) Young tableaux
and multiplication of tableaux.  In particular we shall make use of
the row and column bumping algorithms for tableau multiplication.  For
this and more, see for example \cite{fulton:young}.

A {\em tableau diagram\/} for a set of rank conditions is a filling of
all the boxes in the corresponding rectangle diagram with integers,
such that each rectangle $R_{ij}$ becomes a tableau $T_{ij}$.
Furthermore, it is required that the entries of each tableau $T_{ij}$
are strictly larger than the entries in tableaux above $T_{ij}$ in the
diagram, within 45 degree angles.  These are the tableaux $T_{kl}$
with $i \leq k < l \leq j$ and $(k,l) \neq (i,j)$.

A {\em factor sequence\/} for a tableau diagram with $n$ rows is a
sequence of tableaux $(W_1, \dots, W_n)$, which is obtained as
follows: If $n = 1$ then the only factor sequence is the sequence
$(T_{01})$ containing the only tableau in the diagram.  When $n \geq
2$, a factor sequence is obtained by first constructing a factor
sequence $(U_1, \dots, U_{n-1})$ for the bottom $n-1$ rows of the
tableau diagram, and choosing arbitrary factorizations of the tableaux
in this sequence:
\[ U_i = P_i \cdot Q_i \,. \]
Then the sequence 
\[ (W_1,\dots,W_n) = (T_{01} \cdot P_1 \,,\, 
   Q_1 \cdot T_{12} \cdot P_2 \,, \dots ,\,
   Q_{n-1} \cdot T_{n-1,n})
\]
is the factor sequence for the whole tableau diagram.  The conjecture
from \cite{buch.fulton:chern}, which is the theme of this paper, can
now be stated as follows:

\begin{conj}
\label{conj_orig}
The coefficient $c_\mu(r)$ is equal to the number of different factor
sequences $(W_1, \dots, W_n)$ for any fixed tableau diagram for the
rank conditions $r$, such that $W_i$ has shape $\mu_i$ for each $i$.
\end{conj}

This conjecture first of all implies that the coefficients $c_\mu(r)$
are non-negative and that they are independent of the side lengths of
empty rectangles in the rectangle diagram.  In addition it implies
that the number of factor sequences does only depend on the rectangle
diagram and not on the choice of a filling of its boxes with integers.

\begin{example}
Suppose we are given a sequence of four vector bundles and the
following rank conditions:
\[ \begin{matrix}
E_0 & \to & E_1 & \to & E_2 & \to & E_3
\vspace{0.1cm} \\
1 && 4 && 3 && 3 \\
& 1 && 2 && 2 \\
&& 1 && 1 \\
&&& 0
\end{matrix} \]
These rank conditions then give the following rectangle diagram:
\[ \picT{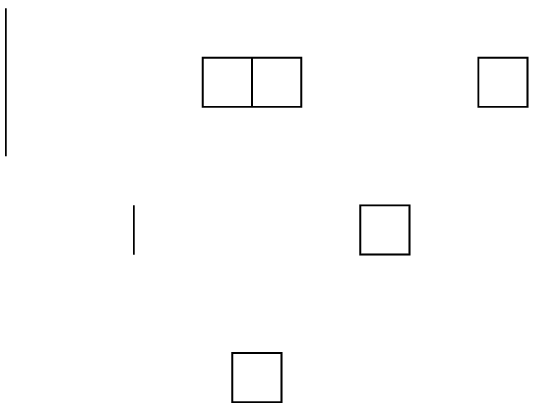} \]
From the bottom row of this diagram we get
\[ P_{\Bar{\Bar r}} = s_{\PART{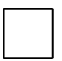}} \,. \]
Then using the algorithm we obtain
\[ P_{\Bar r} = 
   s_{\PART{1}} \otimes s_{\PART{1}} ~+~ 
   1 \otimes s_{\PART{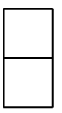}} 
\]
and
\[ \begin{split}
P_r =~ &
s_{\PART{1}} \otimes s_{\PART{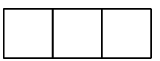}} \otimes s_{\PART{1}} ~+~
s_{\PART{1}} \otimes s_{\PART{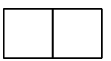}} \otimes s_{\PART{11}} ~+~
1 \otimes s_{\PART{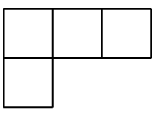}} \otimes s_{\PART{1}} ~+ \\
& 1 \otimes s_{\PART{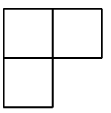}} \otimes s_{\PART{11}} ~+~
1 \otimes s_{\PART{3}} \otimes s_{\PART{11}} ~+~
1 \otimes s_{\PART{2}} \otimes s_{\PART{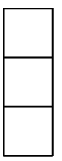}} \,.
\end{split} \]
Thus the formula for the cohomology class of $\Omega_r(E_\bull)$ has
six terms.  Now, one possible tableau diagram for the given rank
conditions is the following:
\[ \picT{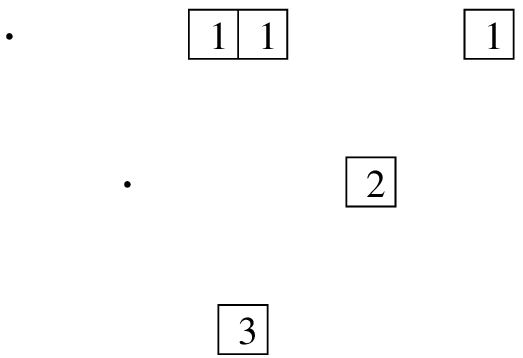} \]
This diagram has the following six factor sequences:
\[ \begin{split}
(\TABLr{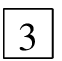}{-3}, \TABLr{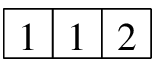}{-3}, \TABLr{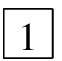}{-3})\,,\,
(\TABLr{3}{-3}, \TABLr{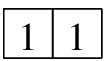}{-3}, \TABLr{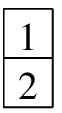}{-3})\,,\,
(\emptyset, \TABLr{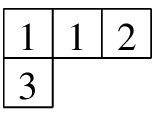}{-3}, \TABLr{1}{-3})\,, \\
(\emptyset, \TABLr{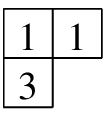}{-3}, \TABLr{21}{-3})\,,\,
(\emptyset, \TABLr{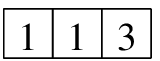}{-3}, \TABLr{21}{-3})\,,\,
(\emptyset, \TABLr{11}{-3}, \TABLr{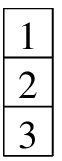}{-3})\,.
\end{split} \]
\end{example}

Since only the rectangle diagram matters for the formula, we will
often depict a rank diagram simply as a triangle of dots in place of a
triangle of numbers.  This is especially convenient when working with
paths through the rank diagram, which we shall do shortly.  Such a
diagram will often be decorated with the rectangles from the
rectangle diagram, or by the tableaux from a tableau diagram.  When
this is done, each rectangle or tableau is put in the middle of the
triangle of dots representing the numbers that produced the rectangle.
In this way the rank conditions used in the above example would be
represented by the diagram:
\[ \picC{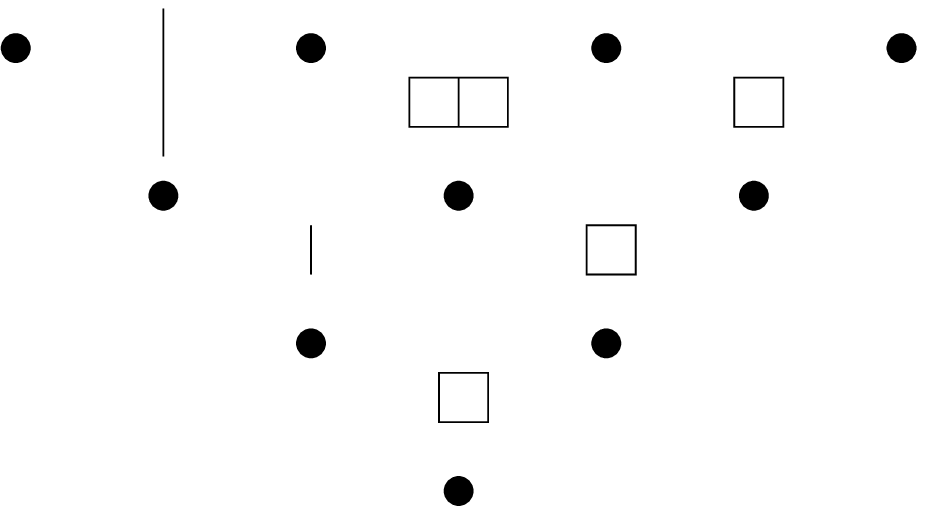} \]

We will now introduce a generalization of the formula $P_r$.  Define a
{\em path\/} through the rank diagram to be a union of line segments
between neighboring rank conditions, which form a continuous path from
$r_{00}$ to $r_{nn}$ such that any vertical line intersects this path
at most once.
\[ \picC{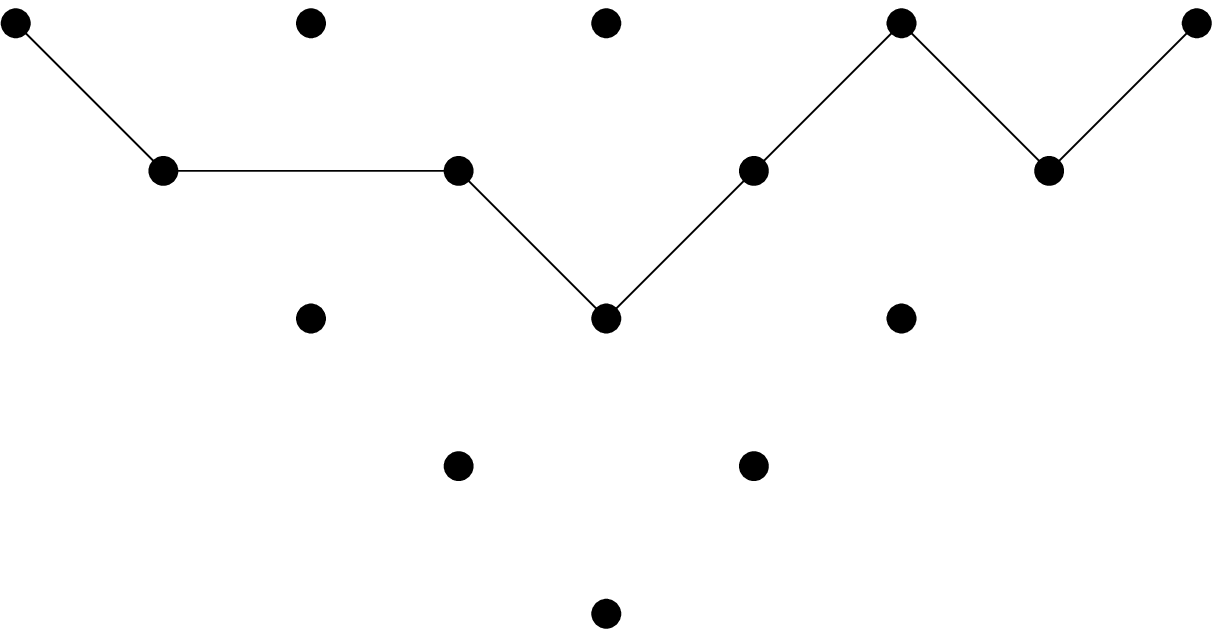} \]
The length of a path is the number of contained line segments (which
is between $n$ and $2n$).  Given a path $\gamma$ of length $\ell$, we
will define an element $P_\gamma \in \Lambda^{\otimes \ell}$.  It is
convenient to identify the basis elements of $\Lambda^{\otimes \ell}$
with labelings of the line segments of $\gamma$ with partitions.  More
generally, if $I_1, \dots, I_\ell$ are sequences of integers, we will
identify the labeling of the line segments in $\gamma$ by these
sequences, left to right, with the element $s_{I_1} \otimes \dots
\otimes s_{I_\ell} \in \Lambda^{\otimes \ell}$.  All basis elements
occurring in $P_\gamma$ will label line segments on the side of the
rank diagram with the empty partition.  If $\gamma$ is the highest
path, going horizontally from $r_{00}$ to $r_{nn}$, then $P_\gamma$ is
equal to $P_r$.

We define $P_\gamma$ inductively as follows.  If $\gamma$ is the
lowest possible path, going from $r_{00}$ to $r_{0n}$ to $r_{nn}$,
then we set $P_\gamma = 1 \otimes 1 \otimes \dots \otimes 1 \in
\Lambda^{\otimes 2n}$.  In other words $P_\gamma$ is equal to the
single basis element which assigns the empty partition to each line
segment.  If $\gamma$ is any other path, then we can find a path
$\gamma'$ which is equal to $\gamma$, except it goes lower at one
place, in one of the following ways:
\begin{center}
\begin{tabular}{ccccc}
& \hspace{.2cm} & $\gamma$ & \hspace{.2cm} & $\gamma'$ 
\vspace{.2cm} \\
Case 1: && \raisebox{-11pt}{\picC{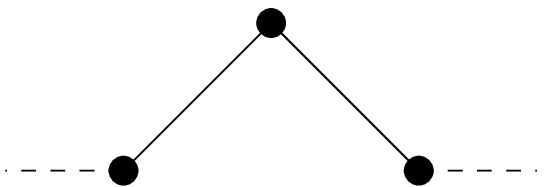}} 
        && \raisebox{-11pt}{\picC{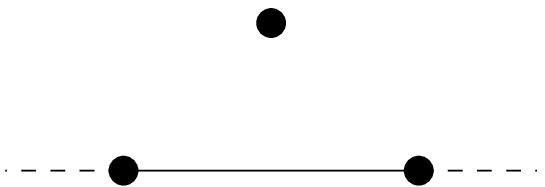}} 
\vspace{.5cm} \\
Case 2: && \raisebox{-11pt}{\picC{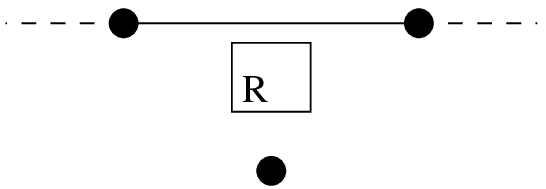}} 
        && \raisebox{-11pt}{\picC{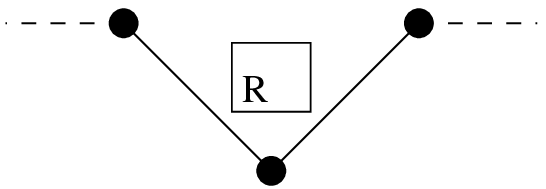}}
\vspace{.2cm}
\end{tabular}
\end{center}
By induction we may assume that $P_{\gamma'}$ is well defined.

If we are in Case 1 we now obtain $P_\gamma$ from $P_{\gamma'}$ by
replacing each basis element
\[ \picC{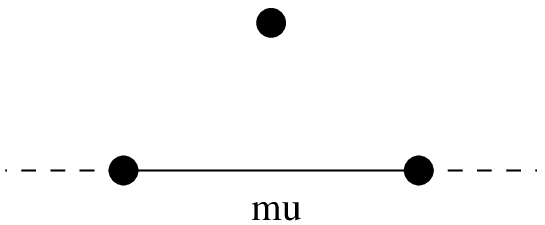} \]
occurring in $P_{\gamma'}$ with the sum
\[ \sum_{\sigma,\tau} c^\mu_{\sigma \tau} 
   \left( \raisebox{-10pt}{\picC{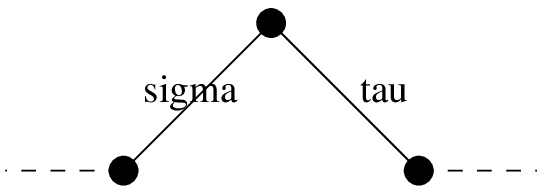}} \right) \,. 
\]
For Case 2, let $R$ be the rectangle associated to the triangle where
$\gamma$ and $\gamma'$ differ.  Then $P_\gamma$ is obtained from
$P_{\gamma'}$ by replacing each basis element
\[ \picC{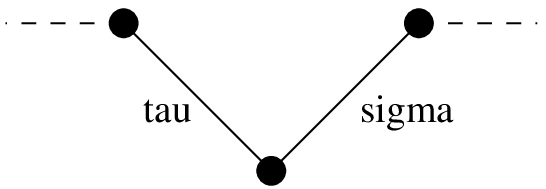} \]
occurring in $P_{\gamma'}$ with zero if $\sigma$ has more rows than
$R$, and otherwise with the element:
\[ \picC{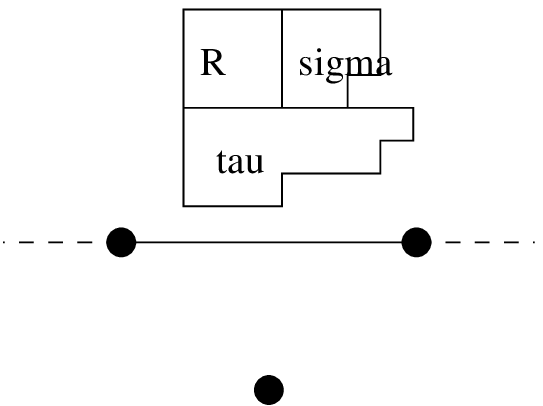} \]
An easy induction shows that this definition is independent of the
choice of $\gamma'$.  The element $P_\gamma$ has geometric meaning
similar to that of $P_r$.  It describes the cohomology class of a
degeneracy locus $\Omega_r(\gamma)$ defined in
\cite{buch.fulton:chern}.

If we are given a tableau diagram, the notion of a factor sequence can
also be extended to paths.  Any factor sequence for a path $\gamma$
will contain one tableau for each line segment in $\gamma$.  As with
elements of $\Lambda^{\otimes \ell}$, we will often regard such a
sequence as a labeling of the line segments in $\gamma$ with tableaux.

If $\gamma$ is the lowest path from $r_{00}$ to $r_{0n}$ to $r_{nn}$
then the only factor sequence is the sequence $(\emptyset, \dots,
\emptyset)$ which assigns the empty tableau to each line segment.
Otherwise we can find a lower path $\gamma'$ as in Case 1 or Case 2
above.  In order to obtain a factor sequence for $\gamma$ we must
first construct one for $\gamma'$.

If we are in Case 1, let $(\dots, W, \dots)$ be a factor sequence for
$\gamma'$ such that $W$ is the label of the displayed line segment,
and let $W = P \cdot Q$ be an arbitrary factorization of $W$.  Then
the sequence $(\dots, P, Q, \dots)$ is a factor sequence for $\gamma$.
For Case 2, let $T$ be the tableau corresponding to the rectangle $R$.
If $(\dots, Q, P, \dots)$ is a factor sequence for $\gamma'$ with $Q$
and $P$ the tableaux assigned to the displayed line segments, then
$(\dots, Q \cdot T \cdot P, \dots)$ is a factor sequence for $\gamma$.

Finally we define coefficients $c_\mu(\gamma) \in \Z$ by the
expression
\[ P_\gamma = \sum_\mu c_\mu(\gamma) \, 
   s_{\mu_1} \otimes \dots \otimes s_{\mu_\ell}  
   \in \Lambda^{\otimes \ell}
\]
where $\ell$ is the length of $\gamma$.  \refconj{conj_orig} then has
the following generalization:

\begin{conj1a}
\label{conj_path}
The coefficient $c_\mu(\gamma)$ is equal to the number of different 
factor sequence $(W_1, \dots, W_\ell)$ for the path $\gamma$, such
that $W_i$ has shape $\mu_i$ for each $i$.
\end{conj1a}


\section{A criterion for factor sequences}
\label{sec_criterion}

In this section we will prove a simple criterion for recognizing
factor sequences.  As in the previous section we will start by
discussing ordinary factor sequences.

Let $\{T_{ij}\}$ be a tableau diagram and let $(W_1, \dots, W_n)$ be a
sequence of tableaux.  At first glance it would appear that to check
if this sequence is a factor sequence, we would have to find all
factor sequences $(U_1, \dots, U_{n-1})$ for the bottom $n-1$ rows of
the tableau diagram, as well as all factorizations $U_i = P_i \cdot
Q_i$, to see if our sequence $(W_1, \dots, W_n)$ is obtained from any
of these, i.e.\ $W_i = Q_{i-1} \cdot T_{i-1,i} \cdot P_i$ for all $i$.
Equivalently we could find all factorizations of each $W_i$ into three
factors $W_i = Q_{i-1} \cdot T_{i-1,i} \cdot P_i$ (with $Q_0 = P_n =
\emptyset$), and check if $(P_1 \cdot Q_1, \dots, P_{n-1} \cdot
Q_{n-1})$ is a factor sequence for any of these choices.  The
criterion for factor sequences allows us to check this for just one
factorization of each $W_i$.

Notice that if the sequence $(W_1, \dots, W_n)$ is a factor sequence,
obtained from an inductive factor sequence $(U_1,\dots, U_{n-1})$ as
above, then the conditions on the filling of a tableau diagram imply
that the entries of each tableau $T_{i-1,i}$ are strictly smaller than
the entries of $Q_{i-1}$ and $P_i$.  This implies that $W_i = Q_{i-1}
\cdot T_{i-1,i} \cdot P_i$ contains the rectangular tableau
$T_{i-1,i}$ in its upper-left corner.
\[ W_i = \raisebox{-20pt}{\picC{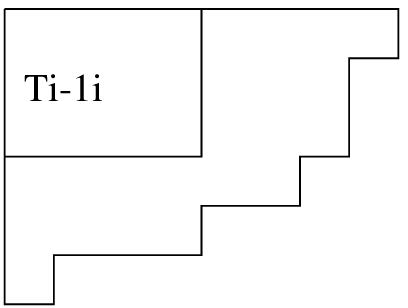}} \]
We shall therefore investigate ways to factor a tableau into three
pieces, one of which is a contained rectangular tableau.

A quick way to factor any tableau is by cutting it along a horizontal
or vertical line.  Let $T$ be a tableau and $a \geq 0$ an integer.
Let $U$ the top $a$ rows of $T$, and $D$ the rest of $T$.  Then $T = D
\cdot U$.  We will call this factorization the {\em horizontal cut\/}
through $T$ after the $a^\Th$ row.  Vertical cuts are defined
similarly.
\begin{align*}
T = \raisebox{-19pt}{\picC{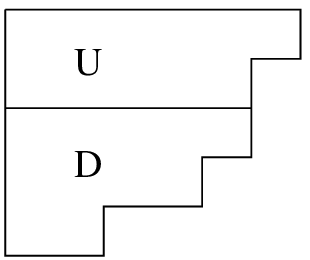}} &= D \cdot U &
T = \raisebox{-19pt}{\picC{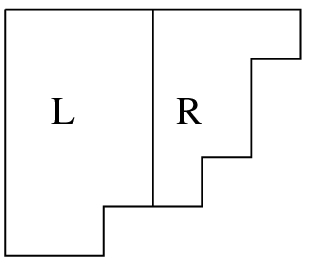}} &= L \cdot R
\end{align*}

\begin{lemma}
\label{lemma_cut}
Let $T = P \cdot Q$ be any factorization of $T$ and let $a$ be the
number of rows in $Q$.  The following are equivalent:
\begin{romenum}
\item $T = P \cdot Q$ is a horizontal cut.
\item The $i^\Th$ row of $T$ has the same number of boxes as the
$i^\Th$ row of $Q$ for $1 \leq i \leq a$.
\item Whenever the top row of $P$ has a box in column $j \geq 1$, the
$a^\Th$ row of $Q$ has a strictly smaller box in this column (unless
$a = 0$).
\end{romenum}
Similarly, if $P$ has $b$ columns, then $T = P \cdot Q$ is a vertical
cut iff the first $b$ columns of $T$ and $P$ have the same heights,
iff the boxes in the last column of $P$ are smaller than or equal to
the boxes in similar positions in the first column of $Q$.
\end{lemma}
\begin{proof}
It is clear that (i) implies (ii) and (iii).  If (iii) is true then
$P$ and $Q$ fit together to form a tableau with $Q$ in the top $a$
rows and $P$ below.  By taking a horizontal cut through this tableau,
we see that it must be the product of $P$ and $Q$.  But then it is
equal to $T$ and (i) follows.  Finally, suppose (ii) is true.  When
the boxes of $P$ are column bumped into $Q$ to form the product $T$,
all of these boxes must then stay below the $a^\Th$ row.  This process
therefore reconstructs $P$ below $Q$ and (i) follows.  The statements
about vertical cuts are proved similarly.
\end{proof}

Now let $W$ be any tableau whose shape contains a rectangle $(b)^a$
with $a$ rows and $b$ columns.  We define the {\em canonical
factorization\/} of $W$ with respect to the rectangle $(b)^a$ to be
the one obtained by first taking a horizontal cut through $W$ after
the $a^\Th$ row, and then a vertical cut through the top part of $W$
after the $b^\Th$ column.
\[ W = \raisebox{-33pt}{\picC{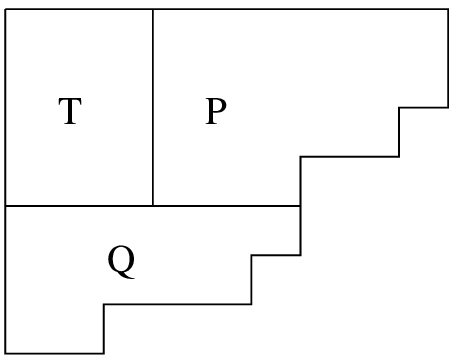}} = Q \cdot T \cdot P \]
Note that this definition depends on $a$, even when $b$ is zero and
the rectangle $(b)^a$ is empty.  When the product of three tableau
$Q$, $T$, and $P$ looks like in this picture, we shall say that the
pair of tableaux $(Q, P)$ {\em fits around\/} the rectangular tableau
$T$.

More generally, let $Q_0$ be the part of $W$ below $T$, $P_0$ the part
of $W$ to the right of $T$, and let $Z$ be the remaining part between
$Q_0$ and $P_0$.
\[ W = \raisebox{-33pt}{\picC{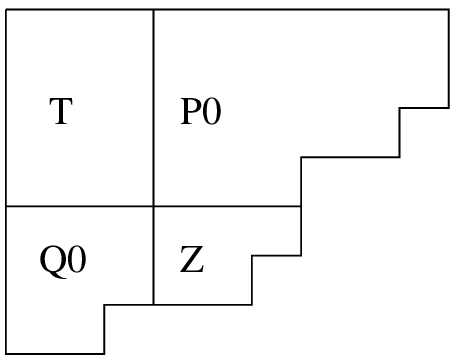}} \]
We define a {\em simple factorization\/} of $W$ with respect to the
rectangle $(b)^a$ to be any factorization $W = Q \cdot T \cdot P$,
such that $Q = Q_0 \cdot \Tilde Q$ and $P = \Tilde P \cdot P_0$ for
some factorization $Z = \Tilde Q \cdot \Tilde P$.

Note that if $Z = \Tilde Q \cdot \Tilde P$ is any factorization of $Z$
and if we put $Q = Q_0 \cdot \Tilde Q$ and $P = \Tilde P \cdot P_0$,
then $Q \cdot T \cdot P = W$.  This follows because $P = \Tilde P
\cdot P_0$ must be a horizontal cut through $P$, and therefore $T
\cdot P = \Tilde P \cdot T \cdot P_0$.  In fact, given arbitrary
tableaux $\Tilde Q$ and $\Tilde P$ one can show that $Q \cdot T \cdot
P = W$ if and only if $\Tilde Q \cdot \Tilde P = Z$, but we shall not
need this here.

We are now ready to formulate the criterion for factor sequences.  Let
$\{R_{ij}\}$ be the rectangles corresponding to the tableau diagram
$\{T_{ij}\}$.  If $(W_1, \dots, W_n)$ is a factor sequence, a simple
factorization of any $W_i$ will always be with respect to the relevant
rectangle $R_{i-1,i}$ from the rectangle diagram.

\begin{thm}
\label{thm_crit}
Let $(W_1, \dots, W_n)$ be a sequence of tableaux such that each $W_i$
contains $T_{i-1,i}$ in its upper-left corner.  Let $W_i = Q_{i-1}
\cdot T_{i-1,i} \cdot P_i$ be any simple factorization of $W_i$ with
respect to the rectangle $R_{i-1,i}$.  Then $(W_1, \dots, W_n)$ is a
factor sequence if and only if $Q_0$ and $P_n$ are empty tableaux and
$(P_1 \cdot Q_1, \dots, P_{n-1} \cdot Q_{n-1})$ is a factor sequence
for the bottom $n-1$ rows of the tableau diagram $\{T_{ij}\}$.
\end{thm}

We shall derive this result from \refprop{prop_pathcrit} below.  Since
this criterion can be applied recursively to the sequence $(P_1 \cdot
Q_1, \dots, P_{n-1} \cdot Q_{n-1})$, it gives an easy algorithm to
determine if a sequence $(W_1, \dots, W_n)$ is a factor sequence.
Note that the easiest way to produce the simple factorizations is to
take the canonical factorization of each $W_i$.  When this choice is
made, the work required in the algorithm essentially consists of
$n(n-1)/2$ tableau multiplications.  Note also that this criterion
makes use of the height of any empty rectangles in the rectangle
diagram.

For proving this criterion we need some definitions.  Let $T$ be a
tableau whose shape is the rectangle $(b)^a$ with $a$ rows and $b$
columns.  We will consider pairs of tableaux $(X,Y)$ such that all
entries in $X$ and $Y$ are strictly larger than the entries of $T$.
For such a pair, let $X = X_0 \cdot \Tilde X$ be the
vertical cut through $X$ after the $b^\Th$ column, and let $Y = \Tilde
Y \cdot Y_0$ be the horizontal cut after row $a$.
\[ \picC{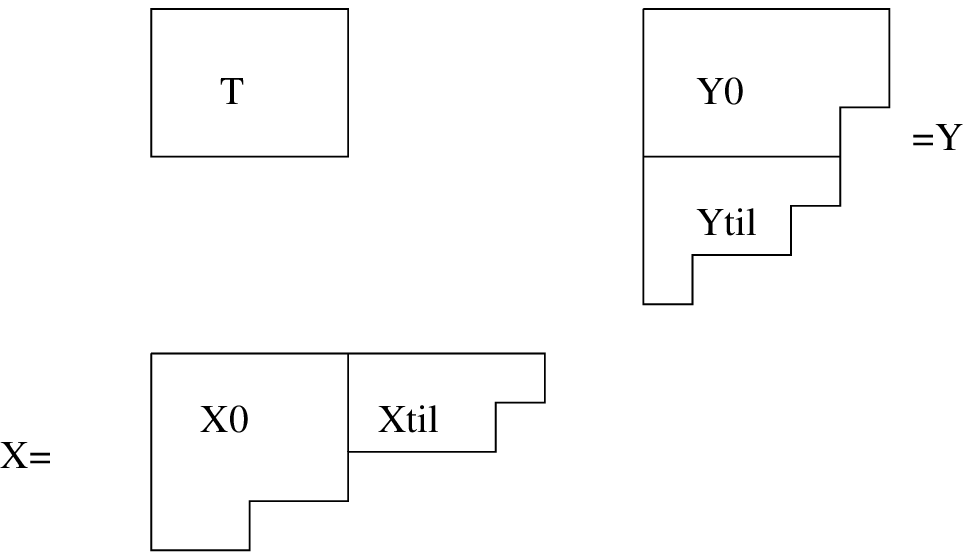} \]
If $(X', Y')$ is another pair of tableaux, we will write $(X, Y)
\models (X', Y')$ if either
\begin{enumerate}
\item for some factorization $\Tilde X = M \cdot N$ we have $X' = X_0
\cdot M$ and $Y' = N \cdot Y$, or
\item for some factorization $\Tilde Y = M \cdot N$ we have $X' = X
\cdot M$ and $Y' = N \cdot Y_0$.
\end{enumerate}
Note that this implies that $X' \cdot T \cdot Y' = X \cdot T \cdot Y$.
In in the first case this follows because $X \cdot T = X_0 \cdot T
\cdot \Tilde X$ and $X' \cdot T = X_0 \cdot T \cdot M$, and the second
case is similar.  We will let $\to$ denote the transitive closure of
the relation $\models$.  This notation depends on the choice of $T$,
as well as the numbers $a$ and $b$ if $T$ is empty.

\begin{lemma}
\label{lemma_tosimple}
Let $W$ be a tableau containing $T$ in its upper-left corner.  Suppose
that the entries of $T$ are smaller than all other entries in $W$.  If
$W = Q \cdot T \cdot P$ is a simple factorization of $W$ with respect
to the rectangle $(b)^a$, and $W = X \cdot T \cdot Y$ is any
factorization, then $(X, Y) \to (Q, P)$.
\end{lemma}
\begin{proof}
Let $X = X_0 \cdot \Tilde X$ be the vertical cut through $X$ after
column $b$, and put $Y' = \Tilde X \cdot Y$.  Then let $Y' = \Tilde Y'
\cdot Y'_0$ be the horizontal cut through $Y'$ after row $a$, and put
$X'' = X_0 \cdot \Tilde Y'$.

We claim that the pair $(X'', Y'_0)$ fits around $T$.  Using
\reflemma{lemma_cut} and that the entries of $T$ are smaller than all
other entries, it is enough to prove that the $b+j^\Th$ entry in the
top row of $X''$ is strictly larger than the $j^\Th$ entry in the
bottom row of $Y'_0$.  This will follow if the $b+j^\Th$ entry in the
top row of $X''$ is larger than or equal to the $j^\Th$ entry in the
top row of $\Tilde Y'$.  Since $X'' = X_0 \cdot \Tilde Y'$ and $X_0$
has at most $b$ columns, this follows from an easy induction on the
number of rows of $\Tilde Y'$.

It follows from the claim that $W = X'' \cdot T \cdot Y'_0$ is the
canonical factorization of $W$, and therefore we have $(X, Y) \models
(X_0, Y') \models (X'', Y'_0) \models (Q,P)$ as required.
\end{proof}

Notice that if $W = X \cdot T \cdot Y$ is a simple factorization and
$(X,Y) \models (X',Y')$, then $W = X' \cdot T \cdot Y'$ must also be a
simple factorization.  It follows that \reflemma{lemma_tosimple} would
be false without the requirement that $W = Q \cdot T \cdot P$ is
simple.

\begin{lemma}
\label{lemma_top}
Let $a \geq 0$ be an integer, and let $Y$ and $S$ be tableaux with
product $A = Y \cdot S$.  Let $A = \Tilde A \cdot A_0$ and $Y = \Tilde
Y \cdot Y_0$ be the horizontal cuts through $A$ and $Y$ after row $a$,
and let $\Tilde Y = M \cdot N$ be any factorization.  Then $N \cdot
Y_0 \cdot S = \Tilde A' \cdot A_0$ for some tableau $\Tilde A'$, and
$M \cdot \Tilde A' = \Tilde A$.
\[ Y = \raisebox{-20pt}{\picC{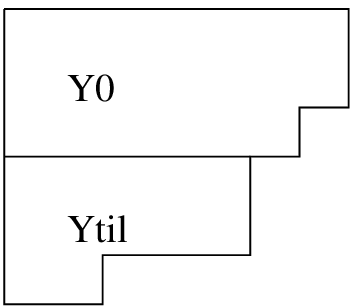}} \;\; ; \;\;\;\; 
   A = Y \cdot S = \raisebox{-20pt}{\picC{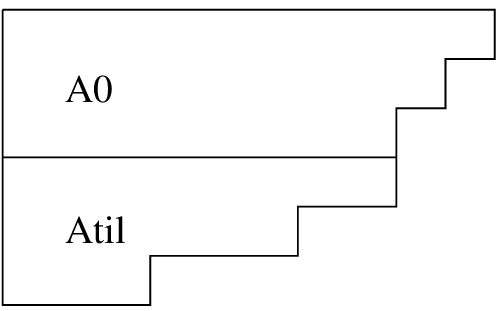}}
\]
\end{lemma}
\begin{proof}
The first statement follows from the observation that the bottom rows
of $Y$ can't influence the top part of $Y \cdot S$, which is a
consequence of the row bumping algorithm.  \reflemma{lemma_cut} then
shows that the factorization $A = (M \cdot \Tilde A') \cdot A_0$ is a
horizontal cut, so $M \cdot \Tilde A' = \Tilde A$ as required.
\end{proof}

\begin{lemma}
\label{lemma_down}
Let $\gamma$ be a path through the rank diagram, and let $(\dots, A, B
\cdot C, \dots)$ be a factor sequence for $\gamma$ such that the
product $B \cdot C$ is the label of a down-going line segment.  Then
$(\dots, A \cdot B, C, \dots)$ is also a factor sequence for $\gamma$.
\[ \picC{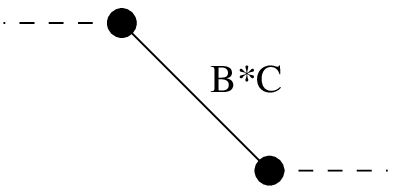} \]
\end{lemma}
\begin{proof}
We will first consider the case where the line segment corresponding
to $A$ goes up.  Let $\gamma'$ be the path under $\gamma$ that cuts
short this line segment and its successor.
\[ \picC{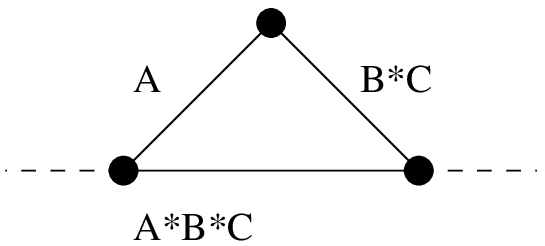} \]
Then by definition $(\dots, A \cdot B \cdot C, \dots)$ is a factor
sequence for $\gamma'$, which means that $(\dots, A \cdot B, C,
\dots)$ is a factor sequence for $\gamma$.  In general $\gamma$ lies
over a path like the one above, and the general case follows from
this.
\end{proof}

Similarly one can prove that if $(\dots, A \cdot B, C, \dots)$ is a
factor sequence for a path, such that $A \cdot B$ is the label of an
up-going line segment, then $(\dots, A, B \cdot C, \dots)$ is also a
factor sequence for this path.

\begin{prop}
\label{prop_pathcrit}
Let $\gamma$ and $\gamma'$ be paths related as in Case 2 of Section
\ref{sec_algo}, and let $(\dots, W, \dots)$ be a factor sequence for
$\gamma$ such that $W$ is the label of the displayed horizontal line
segment.
\[ \picC{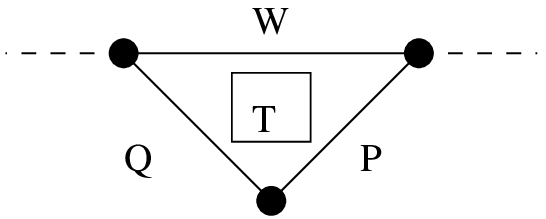} \]
If $W = Q \cdot T \cdot P$ is any simple factorization of $W$, then
$(\dots, Q, P, \dots)$ is a factor sequence for $\gamma'$.
\end{prop}
\begin{proof}
Since $(\dots, W, \dots)$ is a factor sequence for $\gamma$, there
exists a factorization $W = X \cdot T \cdot Y$ such that $(\dots, X,
Y, \dots)$ is a factor sequence for $\gamma'$.  By
\reflemma{lemma_tosimple} we have $(X,Y) \to (Q,P)$.  It is therefore
enough to show that if $(X,Y) \models (X',Y')$ then $(\dots, X', Y',
\dots)$ is a factor sequence for $\gamma'$.

Let $a$ be the number of rows in (the rectangle corresponding to) $T$,
and let $Y = \Tilde Y \cdot Y_0$ be the horizontal cut through $Y$
after the $a^\Th$ row.  We will do the case where a factor of $\Tilde
Y$ is moved to $X$, the other case is proved using a symmetric
argument.  We then have a factorization $\Tilde Y = M \cdot N$ such
that $X' = X \cdot M$ and $Y' = N \cdot Y_0$.  We can assume that the
paths $\gamma$ and $\gamma'$ go down after they meet, and that the
original factor sequence for $\gamma$ is $(\dots, W, S, \dots)$.
\pathskip
\[ \picC{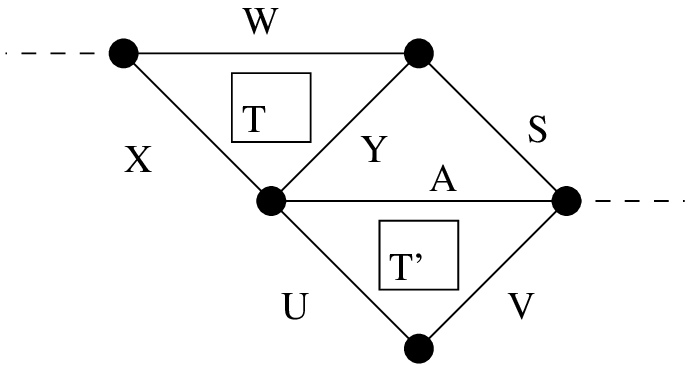} \]
Put $A = Y \cdot S$.  Then $(\dots, X, A, \dots)$ is a factor sequence
for the path with these labels in the picture.  Now let $T'$ be the
rectangular tableau associated to the lower triangle, and let $A = U
\cdot T' \cdot V$ be the canonical factorization of $A$.  Since this
is a simple factorization we may assume by induction that $(\dots, X,
U, V, \dots)$ is a factor sequence.  Using \reflemma{lemma_top} we
deduce that $N \cdot Y_0 \cdot S = U' \cdot T' \cdot V$ for some
tableau $U'$, such that $M \cdot U' = U$.  Since $(\dots, X, M \cdot
U', V, \dots)$ is a factor sequence, so is $(\dots, X \cdot M, U', V,
\dots)$ by \reflemma{lemma_down}.  This means that $(\dots, X \cdot M,
U' \cdot T' \cdot V, \dots) = (\dots, X', Y' \cdot S, \dots)$ is a
factor sequence, which in turn implies that $(\dots, X', Y', S,
\dots)$ is a factor sequence for $\gamma'$ as required.
\end{proof}

The proof of \refprop{prop_pathcrit} also gives the following:

\begin{cor}
Let $(\dots, X, Y, \dots)$ be a factor sequence for the path $\gamma'$
in the proposition.  If $(X, Y) \to (X', Y')$ then $(\dots, X', Y',
\dots)$ is also a factor sequence for $\gamma'$.
\end{cor}
\begin{proof}[Proof of \refthm{thm_crit}]
The ``if'' implication follows from the definition.  If the sequence
$(W_1, \dots, W_n)$ is a factor sequence, then $n$ applications of
\refprop{prop_pathcrit} shows that $(Q_0, P_1, Q_1, P_2, \dots,
Q_{n-1}, P_n)$ is a factor sequence for the path with these labels.
\pathskip
\[ \picC{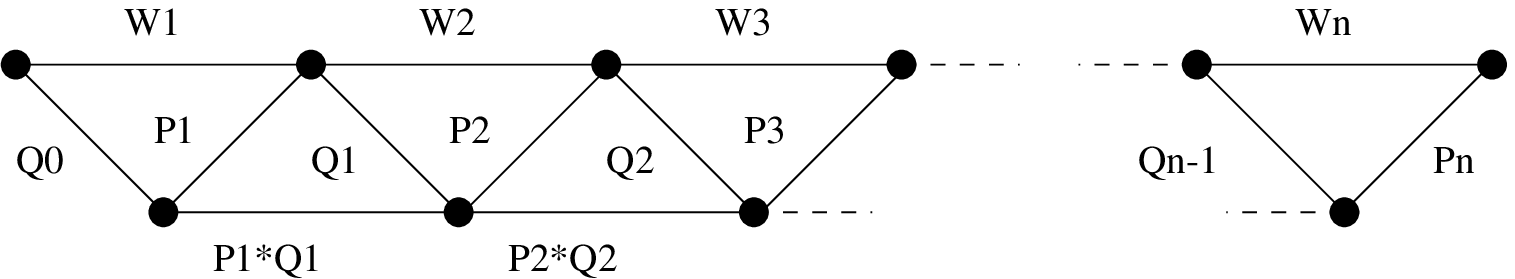} \]
It follows that $Q_0$ and $P_n$ are empty, and $(P_1 \cdot Q_1, \dots,
P_{n-1} \cdot Q_{n-1})$ is a factor sequence for the bottom $n-1$
rows.  This proves ``only if''.
\end{proof}


\section{An involution of Fomin}
\label{sec_fomin}

In this section we will describe a sign-reversing involution on pairs
of tableaux constructed by Sergey Fomin.  The purpose of this
involution is to cancel out the difference between the coefficients
$c_\mu(r)$ produced by the algorithm in Section \ref{sec_algo}, and
their conjectured values.

Fix a positive integer $a$.  If $P$ and $Q$ are tableaux of shapes
$\sigma$ and $\tau$ such that $P$ has at most $a$ rows, we let
$S(\frac{P}{Q})$ denote the symmetric function $s_I \in \Lambda$ where
$I$ is the sequence of integers $I = (\sigma_1, \dots, \sigma_a,
\tau_1, \tau_2, \dots)$.  Let $\PP_a$ be the set of all pairs $(Q,P)$
such that $S(\frac{P}{Q}) \neq 0$ and such that $P$ and $Q$ do not fit
together as a tableau with $P$ in the top $a$ rows and $Q$ below.
This means that the $a^\Th$ row of $P$ must be shorter than the top
row of $Q$, or some box in the top row of $Q$ must be smaller than or
equal to the box in the same position of the $a^\Th$ row of $P$.  For
example, if $a = 2$ the following pairs are in $P_a$:
\[ ( \TABLr{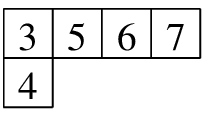}{-5} \,,\, \TABLr{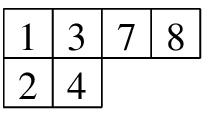}{-5} ) \;\text{ and }\;
   ( \TABLr{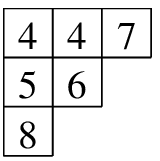}{-5} \,,\, \TABLr{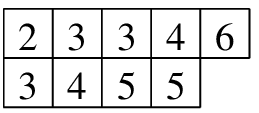}{-5} ) \,.
\]

\begin{lemma}[Fomin's involution]
\label{lemma_fomin}
There exists an involution of $\PP_a$ with the property that if
$(Q, P)$ is mapped to $(Q', P')$ then
\begin{romenum}
\item $Q' \cdot P' = Q \cdot P$,
\item $S(\frac{P'}{Q'}) = - S(\frac{P}{Q})$, and
\item the first column of $Q'$ is equal to the first column of $Q$.
\end{romenum}
\end{lemma}

Fomin supplied the proof of this lemma in the form of the beautiful
algorithm described below.  While Fomin's original description uses
path representations of tableaux, we have translated the algorithm
into notation that is closer to the rest of this paper.

We will work with {\em diagrams with weakly increasing rows\/}.  These
will be ``Young diagrams'' for finite sequences of non-negative
integers, where all boxes are filled with integers so that the rows
are weakly increasing.  Empty rows are allowed as in the following
example:
\[ \picT{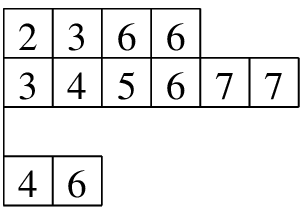} \]
A {\em violation\/} for such a diagram to be a tableau is a box in the
second row or below, such that there is no box directly above it, or
the box directly above it is not strictly smaller.  The above
diagram has 4 violations in its second row and 2 in row four.

If $D$ is a diagram with weakly increasing rows, and if $I$ is the
sequence of row lengths, we put $S(D) = s_I \in \Lambda$.  Let
$\rect(D)$ denote the tableau obtained by multiplying the rows of $D$
together, from bottom to top.  We will identify a pair $(Q,P) \in
\PP_a$ with the diagram $D$ consisting of $P$ in the top $a$ rows and
$Q$ below.  For this diagram we then have $Q \cdot P = \rect(D)$ and
$S(\frac{P}{Q}) = S(D)$.

We will start by taking care of the special case where $a=1$ and both
$P$ and $Q$ have at most one row.  In this case \reflemma{lemma_fomin}
without property (iii) is equivalent to the identity $s_{\ell,k} =
h_\ell h_k - h_{\ell+1} h_{k-1}$ in the plactic monoid, which is a
special case of a result by Lascoux and Sch{\"u}tzenberger
\cite{schutzenberger:correspondance},
\cite{lascoux.schutzenberger:mono}.  The simple proof of this result
given in \cite{fomin.greene:noncommutative} develops techniques which
Fomin used to establish \reflemma{lemma_fomin} in full generality.

\begin{lemma}
\label{lemma_tworows}
Let $D$ be a diagram with two rows and at least one violation in the
second row.  Then there exists a unique diagram $D'$ such that
$\rect(D') = \rect(D)$ and $S(D') = -S(D)$.  Furthermore, $D'$ also
has two rows and at least one violation in the second row.  The
leftmost violations of $D$ and $D'$ appear in the same column and
contain the same number.  The parts of $D$ and $D'$ to the left of
this column agree.
\end{lemma}
\begin{proof}
Let $p$ and $q$ be the lengths of the top and bottom rows of $D$.  The
requirement $S(D') = -S(D)$ then implies that $D'$ must have two rows
with $q-1$ boxes in the top row and $p+1$ in the bottom row.  Now it
follows from the Pieri formula \cite[\S 2.2]{fulton:young} that the
product $\rect(D)$ of the rows in $D$ has at most two rows.
Furthermore, since $D$ contains a violation, the second row of
$\rect(D)$ has at most $q-1$ boxes.  Using the Pieri formula again,
this implies that there is exactly one way to factorize $\rect(D)$
into a row of length $p+1$ times another of length $q-1$.  This
establishes the existence and uniqueness of $D'$.

Explicitly, one may use the inverse row bumping algorithm to obtain
this factorization of $\rect(D)$.  This is done by bumping out a
horizontal strip of $q-1$ boxes which includes all boxes in the second
row, working from right to left.

Let $x$ be the leftmost violation of $D$, where $D$ has the form:
\[ D = \raisebox{-9pt}{\picT{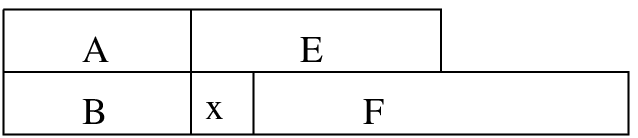}} \,. \]
Suppose the parts $A$ and $B$ each contain $t$ boxes.  Now form the
product $F \cdot E$ and let $c_j$ and $d_j$ be the boxes of this
product as in the picture:
\[ F \cdot E = \raisebox{-9pt}{\picT{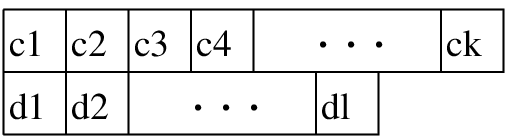}} \,. \]
Since $x$ is a violation in $D$, it must be smaller than all boxes in
$E$ and $F$.  Therefore we have
\[ x \cdot F \cdot E = \raisebox{-9pt}{\picT{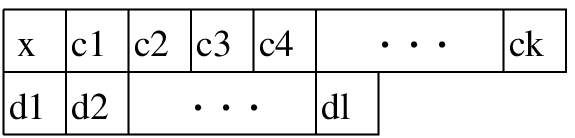}} \,. \]
Now since each $d_j > c_j$ it follows that if a horizontal strip of
length $q-t-1$ is bumped off this tableau, $x$ will remain where it
is.  In other words we can factor $x \cdot F \cdot E$ into $x \cdot F'
\cdot E'$ such that $x \cdot F'$ and $E'$ are rows of lengths $p-t+1$
and $q-t-1$ respectively.  Since the entries of $A$ and $B$ are no
larger than $x$, the products $B \cdot x \cdot F'$ and $A \cdot E'$
are rows of lengths $p+1$ and $q-1$.  But the product of these rows is
$\rect(D)$, so they must be the rows of $D'$ by the uniqueness.  This
proves that $D'$ has the stated properties.
\end{proof}

Notice that the uniqueness also implies that the transformation of
diagrams described in the lemma is inverse to itself, i.e.\ an
involution.

Now suppose $D$ is any diagram with weakly increasing rows.  Then
\reflemma{lemma_tworows} can be applied to any subdiagram of two
consecutive rows, such that the second of these rows contains a
violation.  If this subdiagram is replaced by the new two-row diagram
given by the lemma, we arrive at a diagram $D'$ satisfying $S(D') =
-S(D)$ and $\rect(D') = \rect(D)$.  We will call this an {\em exchange
operation\/} between the two rows of $D$.

We shall need an ordering on the violations in a diagram.  Here the
smallest of two violations is the south-west most one.  If the two
violations are equally far south-west, then the north-west most one is
smaller.  In other words, a violation in row $i$ and column $j$ is
smaller than another in row $i'$ and column $j'$ iff $j-i < j'-i'$,
or $j-i = j'-i'$ and $i < i'$.

Notice that when an exchange operation between two rows is carried
out, violations may appear or disappear in these two rows as well as
in the row below them.  However, the properties given in
\reflemma{lemma_tworows} imply that all of the changed violations will
be larger than the left-most violation in the second of the rows
exchanged.  It follows that the minimal violation in a diagram will
remain constant if any (sequence of) exchange operations is carried
out.  Similarly, all boxes south-west of the minimal violation will
remain fixed.

\newcommand{\DD}{{\mathcal D}_{Q,P}}

\begin{proof}[Proof of \reflemma{lemma_fomin}] 
Given a pair $(Q, P) \in \PP_a$, let $\DD$ be the finite set of all
non-tableau diagrams $D$ with weakly increasing rows, such that
$\rect(D) = Q \cdot P$ and $S(D) = \pm S(\frac{P}{Q})$, and so that
the minimal violation in $D$ is in row $a+1$.  The pair $(Q, P)$ is
then identified with one of the diagrams in this set.  We will
describe an involution of the set $\DD$ and another of the complement
of $\PP_a \cap \DD$ in $\DD$.  The restriction of Fomin's involution
to $\PP_a \cap \DD$ is then obtained by applying the involution
principle of Garsia and Milne \cite{garsia.milne:method} to these
involutions.

The involution of $\DD$ simply consists of doing an exchange operation
between the rows $a$ and $a+1$ of a diagram.  This is possible because
all diagrams are required to have a violation in row $a+1$.

Now note that a diagram $D \in \DD$ is in the complement of $\PP_a
\cap \DD$ if and only if $D$ has a violation outside the $a+1^\St$
row.  We take the involution of $\DD \setminus \PP_a$ to be an
exchange operation between the row of the minimal violation outside
row $a+1$, and the row above this violation.  This is indeed an
involution since the minimal violation outside row $a+1$ stays the
same.

These involutions now combine to give an involution of $\PP_a \cap
\DD$ by the involution principle.  To carry it out, start by forming
the diagram with $P$ in the top $a$ rows and $Q$ below it.  Then do an
exchange operation between row $a$ and row $a+1$.  If all violations
in the resulting diagram are in row $a+1$ we are done.  $P'$ is then
the top $a$ rows of this diagram and $Q'$ is the rest.  Otherwise we
continue by doing an exchange operation between the row of the minimal
violation outside row $a+1$ and the row above it, followed by another
exchange operation between row $a$ and row $a+1$.  We continue in this
way until all violations are in row $a+1$.

Finally, the properties of $P'$ and $Q'$ follow from the properties of
exchange operations.  In particular, the requirement $S(\frac{P'}{Q'})
= -S(\frac{P}{Q})$ follows because we always carry out an odd number
of exchange operations.
\end{proof}

\begin{example}
  The pair $(P, Q) = (\TABLr{43567}{-5} \,,\, \TABLr{241378}{-5})$ in $\PP_2$
  gives the following sequence of exchange operations:
\[ \raisebox{-12pt}{\TABL{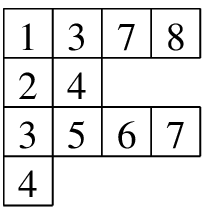}} \;\rightsquigarrow\;
   \raisebox{-12pt}{\TABL{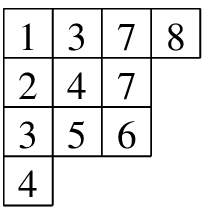}} \;\rightsquigarrow\;
   \raisebox{-12pt}{\TABL{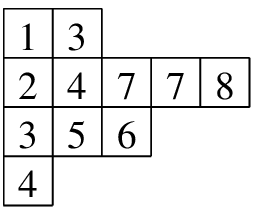}} \;\rightsquigarrow\;
   \raisebox{-12pt}{\TABL{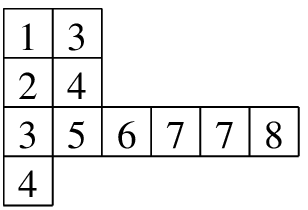}}
\]
This pair therefore corresponds to $(P', Q') = (\TABLr{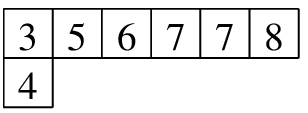}{-5} \,,\,
\TABLr{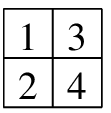}{-5})$ by Fomin's involution.
\end{example}

There are examples of pairs $(Q, P)$ for which the set $\PP_a \cap
\DD$ has more than two diagrams, all with the same first column.  This
means that the involution constructed above is not the only one that
satisfies the conditions of \reflemma{lemma_fomin}.  One way to
produce different involutions is to use another ordering among
violations.  The only property of the order that we have used is that
when an exchange operation is performed, any appearing and
disappearing violations must be larger than the leftmost violation in
the second of the rows being exchanged.  For example, given any
irrational parameter $\xi \in (0,1)$, we obtain a new order by letting
a violation in position $(i,j)$ be smaller than another in position
$(i',j')$ if and only if $j - \xi i < j' - \xi i'$.



\section{The stronger conjecture}
\label{sec_stronger}

In this section we will present a simple conjecture which implies
Conjecture~1A.  Let $\gamma$ be a path through the rank diagram which
at some triangle has an angle pointing down:
\[ \picC{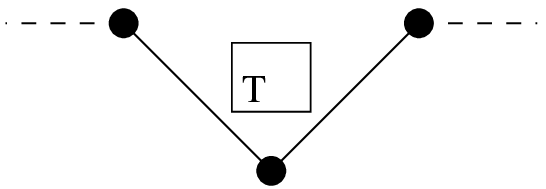} \]
Let $T$ be the rectangular tableau associated to this triangle, and
suppose the corresponding rectangle has $a$ rows and $b$ columns.

If $X$ and $Y$ are tableaux whose entries are strictly larger than the
entries of $T$, and if $Y$ has at most $a$ rows, we will let
\[ \attach{T}{X}{Y} = \raisebox{-25pt}{\picC{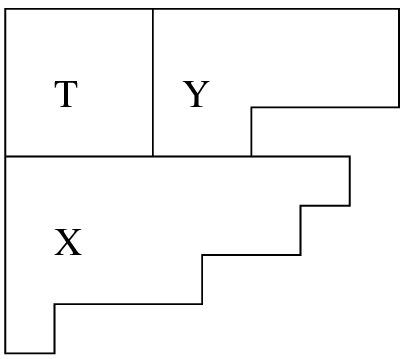}} \]
denote the diagram consisting of $T \cdot Y$ in the top $a$ rows and
$X$ below.  The sequence of row lengths of this diagram then gives an
element $S(\attach{T}{X}{Y})$ in the ring of symmetric functions
$\Lambda$.  Note that $(X, Y)$ fits around $T$ if and only if the
diagram $\attach{T}{X}{Y}$ is a tableau.

Suppose that $(X,Y)$ does not fit around $T$ and $S(\attach{T}{X}{Y})$
is non-zero.  Let $X = X_0 \cdot \Tilde X$ be the vertical cut through
$X$ after the $b^\Th$ column.  Then $(\Tilde X, Y)$ is an element of
the set $\PP_a$ defined in the previous section.  Let $(\Tilde X',
Y')$ be the result of applying Fomin's involution to this pair, and
set $X' = X_0 \cdot \Tilde X'$.  Since the first columns of $\Tilde X$
and $\Tilde X'$ agree, $X'$ consists of $X_0$ with $\Tilde X'$
attached to its right side by \reflemma{lemma_cut}.  It follows that
$S(\attach{T}{X'}{Y'}) = -S(\attach{T}{X}{Y})$.  (Note that one could
also get from $(X, Y)$ to $(X', Y')$ by applying Fomin's involution to
the pair $(X, T \cdot Y)$.)

\begin{conj}
\label{conj_strong}
Let $(\dots, X, Y, \dots)$ be a factor sequence for $\gamma$ with $X$
and $Y$ the labels of the displayed line segments, such that $Y$ has
at most $a$ rows.  Suppose $(X, Y)$ does not fit around $T$ and
$S(\attach{T}{X}{Y}) \neq 0$.  If $X'$ and $Y'$ are obtained from $X$
and $Y$ by applying Fomin's involution as described above, then
$(\dots, X', Y', \dots)$ is also a factor sequence for $\gamma$.
\end{conj}

If we fix the location of the down-pointing angle of $\gamma$ (i.e.\
the location of $T$ in the tableau diagram), then the strongest case
of this conjecture is when the rest of $\gamma$ goes as low as
possible.  If \refconj{conj_strong} is true for all locations of the
down-pointing angle, then the conjectured formula for the coefficients
$c_\mu(\gamma)$ is correct.

\begin{thm}
\label{thm_strong_implies}
Conjecture~1A follows from \refconj{conj_strong}.
\end{thm}
\begin{proof}
If $W_1, \dots, W_\ell$ are diagrams with weakly increasing rows,
e.g.\ tableaux, we will write $S(W_1, \dots, W_\ell) = S(W_1) \otimes
\dots \otimes S(W_\ell) \in \Lambda^{\otimes \ell}$.  With this
notation we must prove that if $\gamma$ is a path of length $\ell$,
then
\begin{equation}
\label{eqn_strong_toprove}
P_\gamma = \sum_{(W_i)} S(W_1, \dots, W_\ell)
\end{equation}
where the sum is over all factor sequences $(W_i)$ for $\gamma$.

Let $\gamma'$ be a path under $\gamma$ as in Case 1 or Case 2 of
Section \ref{sec_algo}.  By induction we can assume that Conjecture~1A
is true for $\gamma'$, i.e.
\begin{equation}
\label{eqn_strong_induct}
P_{\gamma'} = \sum_{(U_i)} S(U_1, \dots, U_{\ell'}) 
\end{equation}
where this sum is over the factor sequences for $\gamma'$.  We must
prove that the right hand side of (\ref{eqn_strong_toprove}) is
obtained by replacing each basis element of (\ref{eqn_strong_induct})
in the way prescribed by the definition of $P_\gamma$.  If we are in
Case 1 then this follows from the Littlewood-Richardson rule \cite[\S
5.1]{fulton:young}: If $U$ is a tableau of shape $\mu$ and $\sigma$
and $\tau$ are partitions, then there are $c^\mu_{\sigma \tau}$ ways
to factor $U$ into a product $U = P \cdot Q$ such that $P$ has shape
$\sigma$ and $Q$ has shape $\tau$.

Assume we are in Case 2.  By induction we then have $P_{\gamma'} =
\sum S(\dots, X, Y, \dots)$ where the sum is over all factor sequences
$(\dots, X, Y, \dots)$ for $\gamma'$; $X$ and $Y$ are the labels of
the two line segments where $\gamma'$ is lower than $\gamma$.  Let $T$
be the rectangular tableau of the corresponding triangle, and let $a$
be the number of rows in its rectangle.  Then by definition we get
\begin{equation}
\label{eqn_strong_extra}
P_\gamma = \sum_{(\dots, X, Y, \dots)} 
S(\dots, \attach{T}{X}{Y}, \dots) 
\end{equation}
where the sum is over all factor sequences $(\dots, X, Y, \dots)$ for
$\gamma'$ such that $Y$ has at most $a$ rows.

Now suppose we have a factor sequence $(\dots, X, Y, \dots)$ such that
the diagram $\attach{T}{X}{Y}$ is a tableau.  Then this tableau must
be the product $X \cdot T \cdot Y$, and so $(\dots, \attach{T}{X}{Y},
\dots)$ is a factor sequence for $\gamma$.  Thus the term $S(\dots,
\attach{T}{X}{Y}, \dots)$ matches one of the terms of
(\ref{eqn_strong_toprove}).  On the other hand it follows from
\refprop{prop_pathcrit} that every term of (\ref{eqn_strong_toprove})
is matched in this way.

We conclude from this that the terms in (\ref{eqn_strong_toprove}) is
the subset of the terms in (\ref{eqn_strong_extra}) which come from
factor sequences such that $(X,Y)$ fits around $T$.  We claim that the
sum of the remaining terms in (\ref{eqn_strong_extra}) is zero.  In
fact, if $(\dots, X, Y, \dots)$ is a factor sequence for $\gamma'$
such that $(X, Y)$ doesn't fit around $T$ and $S(\attach{T}{X}{Y})
\neq 0$, then we may apply Fomin's involution in the way described
above to get tableaux $X'$ and $Y'$.  If \refconj{conj_strong} is
true, then the sequence $(\dots, X', Y', \dots)$ is also a factor
sequence, and since $S(\attach{T}{X'}{Y'}) = - S(\attach{T}{X}{Y})$,
the terms of (\ref{eqn_strong_extra}) given by these two factor
sequences cancel each other out.
\end{proof}

The number of factor sequences for a tableau diagram can be extremely
large.  For this reason it is almost impossible to verify
\refconj{conj_orig} or Conjecture~1A by computing both sides of their
equations.  In contrast, instances of \refconj{conj_strong} can be
tested easily even on large examples.  Given a tableau diagram and a
path, one can generate a factor sequence for this path by choosing
factorizations of tableaux by random.  Then one can apply Fomin's
involution to the sequence, and use the criterion of
\refprop{prop_pathcrit} to check that the result is still a factor
sequence.  Such checks have been carried out repeatedly for each of
500,000 randomly chosen tableau diagrams with up to 10 rows of
tableaux, without finding any violations of \refconj{conj_strong}.
Together with the results in the next section, we consider this to be
convincing evidence for the conjectures.


\section{Proof in a special case}
\label{sec_fourbdl}

In this final section we will show that \refconj{conj_strong} is true
in certain special cases.  These cases will be sufficient to prove the
conjectured formula for $c_\mu(r)$ when all rectangles in and below
the fourth row of the rectangle diagram are empty, and when no two
non-empty rectangles in the third row are neighbors.  This covers all
situations with at most four vector bundles.

Let $\gamma$ be a path through the rank diagram with a down-pointing
angle as in the previous section.  Let $R$ be the rectangle of the
corresponding triangle.
\[ \picC{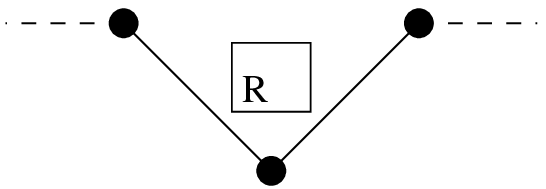} \]
We will describe two cases where \refconj{conj_strong} can be proved.
Both cases require a special configuration of the rectangles
surrounding $R$.  Suppose $R$ is the rectangle $R_{ij}$ in the
rectangle diagram.  We will say that a different rectangle $R' =
R_{kl}$ is {\em below\/} $R$ if $k \leq i < j \leq l$.  $R'$ is {\em
strictly below\/} $R$ if $k < i < j < l$.

\begin{prop}
\label{prop_nobelow}
\refconj{conj_strong} is true for $\gamma$ if all rectangles strictly
below $R$ are empty.
\[ \picC{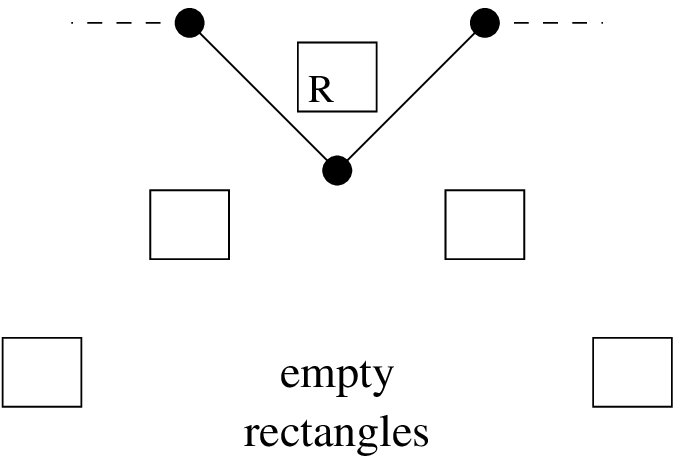} \]
\end{prop}
Note that this covers all rectangles on the left and right sides of
the rectangle diagram.
\begin{proof}
Let $T$ be the tableau corresponding to $R$, and suppose $(\dots, X,
Y, \dots)$ is a factor sequence for $\gamma$.  Since all tableau on
the line going south-west from $T$ in the tableau diagram are narrower
than $T$, it follows that also $X$ has fewer columns than $T$.
Similarly $Y$ has fewer rows than $T$.  But this means that $(X,Y)$
fits around $T$ and the statement of \refconj{conj_strong} is
trivially true.
\end{proof}

In the other situation we shall describe, we allow three non-empty
tableaux below $T$ as shown in the picture.
\vspace{3pt}
\[ \picC{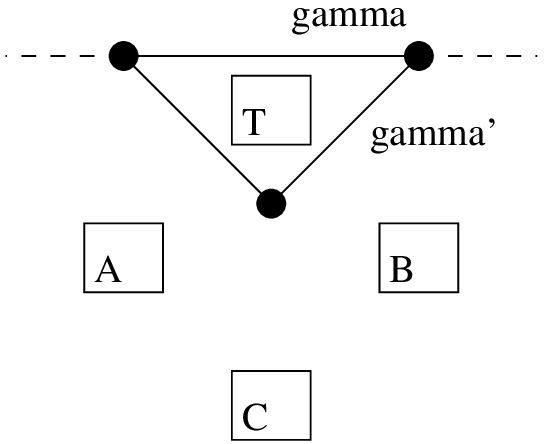} \]
All other tableaux below $T$ are required to be empty.  Let $\gamma$
be the higher and $\gamma'$ the lower of the two paths in the
diagram.

\begin{lemma}
\label{lemma_fourbdl}
Let $(\dots, X, Y, \dots)$ be a labeling of the line segments of
$\gamma'$ with tableaux.  The following are equivalent:
\begin{mathenum}
\item $(\dots, X, Y, \dots)$ is a factor sequence for $\gamma'$.
\item $(\dots, X \cdot T \cdot Y, \dots)$ is a factor sequence for
$\gamma$ and the part of $X$ that is wider than $T$ and the part of
$Y$ that is taller than $T$ have entries only from $C$.
\end{mathenum}
\end{lemma}
\begin{proof}
It is clear that (1) implies (2).  For the other implication, put $W =
X \cdot T \cdot Y$ and let $W = X' \cdot T \cdot Y'$ be the canonical
factorization of $W$.  Then it follows from \refprop{prop_pathcrit}
that $(\dots, X', Y', \dots)$ is a factor sequence for $\gamma'$.
Since $(X, Y) \to (X', Y')$ by \reflemma{lemma_tosimple}, we may
assume that $(X, Y) \models (X', Y')$.

We will handle the case where a factor of the bottom part of $Y$ is
moved to $X$, the other case being symmetric.  This means that for
some tableau $M$ we have $X' = X \cdot M$ and $Y = M \cdot Y'$.  Since
the bottom part of $Y$ has entries only from $C$, this is also true
for $M$.

We may assume that $\gamma$ and $\gamma'$ go down outside the
displayed angle and that our factor sequence is $(\dots, U, X', Y', V,
\dots)$.
\pathskip
\[ \picC{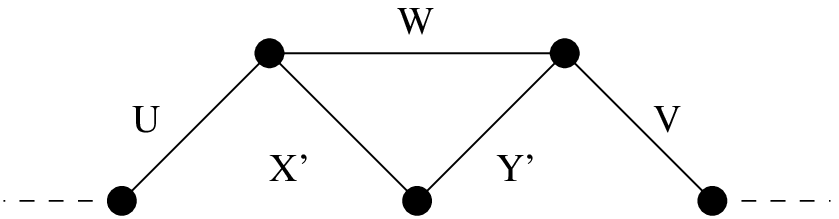} \]
Then by definition there exists a factorization $C = C'_1 \cdot C'_2$
such that $A \cdot C'_1 = U \cdot X'$ and $C'_2 \cdot B = Y' \cdot V$.
Since $U \cdot X \cdot M$ consists of $A$ with $C'_1$ attached on its
right side, and since all entries of $M$ are strictly larger than the
entries of $A$, it follows that $U \cdot X$ consists of $A$ with some
tableau $C_1$ attached on the right side.  Furthermore $C_1 \cdot M =
C'_1$ by \reflemma{lemma_cut}.

Put $C_2 = M \cdot C'_2$.  Then we have $C_1 \cdot C_2 = C$, $A \cdot
C_1 = U \cdot X$, and $C_2 \cdot B = Y \cdot V$.  It follows that
$(\dots, U, X, Y, V, \dots)$ is a factor sequence as required.
\end{proof}

\begin{prop}
\label{prop_onebelow}
\refconj{conj_strong} is true for the path $\gamma'$ in
\reflemma{lemma_fourbdl}.
\end{prop}
\begin{proof}
Let $(\dots, X, Y, \dots)$ be a factor sequence for $\gamma'$ which
satisfies the conditions in \refconj{conj_strong}, and let $X'$ and
$Y'$ be the tableaux obtained from $X$ and $Y$ using Fomin's
involution.  Since the part of $X$ that is wider than $T$ has entries
only from $C$, the same will be true for $X'$ by
\reflemma{lemma_fomin} (iii).  Since $Y'$ has fewer rows than $T$ and
since $(\dots, X' \cdot T \cdot Y', \dots) = (\dots, X \cdot T \cdot
Y, \cdots)$ is a factor sequence for $\gamma$, it follows from
\reflemma{lemma_fourbdl} that $(\dots, X', Y', \dots)$ is a factor
sequence for $\gamma'$.
\end{proof}

\begin{cor}
  \refconj{conj_orig} is true if all rectangles in and below the
  fourth row of the rectangle diagram are empty, and if no two
  non-empty rectangles in the third row are neighbors.
\end{cor}
\begin{proof}
When the rectangle diagram satisfy these properties, then all
instances of \refconj{conj_strong} follow from either
\refprop{prop_nobelow} or \refprop{prop_onebelow}.  The corollary
therefore follows from \refthm{thm_strong_implies}.
\end{proof}

In Section \ref{sec_algo} we defined a rectangle diagram to be
something you get by replacing the small triangles of numbers in a
rank diagram with rectangles.  However, everything we have done is
still true if one defines a rectangle diagram to be any diagram of
rectangles, each given by a number of rows and columns, such that the
number of rows decreases when one moves south-east while the number of
columns decreases when one moves south-west.  This definition is
slightly more general because the side lengths of the rectangles in a
rectangle diagram obtained from rank conditions satisfy certain
relations.  Although we don't know any geometric interpretation of the
more general rectangle diagrams, they seem to be the natural
definition for combinatorial purposes.



\providecommand{\bysame}{\leavevmode\hbox to3em{\hrulefill}\thinspace}


\end{document}